\documentclass[12pt, a4paper, oneside]{article}
\usepackage[top=1in, bottom=1.1in, left=0.9in, right=0.5in]{geometry}
\usepackage[parfill]{parskip}
\usepackage{amsmath}
\usepackage{amsthm}
\usepackage{graphicx}
\usepackage{amssymb}
\usepackage{epstopdf}
\usepackage{setspace}
\usepackage{authblk}
\usepackage{enumerate}
\usepackage{lmodern}
\usepackage[hypcap]{caption}
\usepackage[english]{babel}
\usepackage{fancyhdr}
\addto\captionsenglish{}
\addto\captionsenglish{}
\addto\captionsenglish{}
\addto\captionsenglish{}
\addto\captionsenglish{}

\newlength\longest

\begin{document}
	
	\newtheorem{theorem}{Theorem}[section]
	\newtheorem{lemma}[theorem]{Lemma}
	\newtheorem{corollary}[theorem]{Corollary}
	\newtheorem{proposition}[theorem]{Proposition}
	\newtheorem{conjecture}[theorem]{Conjecture}
	\newtheorem{problem}[theorem]{Problem}
	\newtheorem{claim}[theorem]{Claim}
	\theoremstyle{definition}
	\newtheorem{assumption}[theorem]{Assumption}
	\newtheorem{fact}[theorem]{Fact}
	\newtheorem{remark}[theorem]{Remark}
	\newtheorem{definition}[theorem]{Definition}
	\newtheorem{example}[theorem]{Example}
	\theoremstyle{remark}
	\newtheorem{notation}{Notasi}
	\renewcommand{\thenotation}{}

\title{A Formula for The $g$-Angle between Two Subspaces of a Normed Space}
\author{M. Nur${}^{1}$, H. Gunawan${}^{2}$ and O. Neswan${}^{3}$}
\affil{${}^{1,2,3}$Analysis and Geometry Group, Faculty of Mathematics and\\
	Natural Sciences, Bandung Institute of Technology,\\
	Jl. Ganesha 10, Bandung 40132, Indonesia\\
\bigskip
E-mail: ${}^{1}$nur$_{-}$math@student.itb.ac.id, ${}^{2}$hgunawan@math.itb.ac.id,
${}^{3}$oneswan@math.itb.ac.id}
\date{}

\maketitle
\begin{abstract}
\noindent We develop the notion of $g$-angle between two subspaces of a normed space.
In particular, we discuss the $g$-angle between a $1$-dimensional subspace and a
$t$-dimensional subspace for $t\geq 1$ and the $g$-angle between a $2$-dimensional
subspace and a $t$-dimensional subspace for $t\geq 2$. Moreover, we present an
explicit formula for the $g$-angle between two subspaces of $\ell^{p}$ spaces.

\bigskip{}

\bigskip{}
\end{abstract}

\section{Introduction}

In an inner product space $(X,\left\langle\cdot,\cdot\right\rangle)$, we can calculate
the angles between two vectors and two subspaces. In particular, the angle $\theta=\theta(x,y)$
between two nonzero vectors $x$ and $y$ in $X$ is defined by
$\cos \theta :=\frac{\left\langle x,y\right\rangle}{\|x\| \| y\| }$
where $\| x\|:=\left\langle x,x\right\rangle^\frac{1}{2}$ denotes the induced norm on $X$.
One may observe that the angle $\theta$ in $X$ satisfies the following basic properties
(see \cite{Diminnie}).

\medskip

\begin{enumerate}
\item[(a)] \emph{Parallelism}: $\theta(x,y)=0$ if and only if $x$ and $y$ are of the same
direction; $\theta(x,y)=\pi $ if and only if $x$ and $y$ are of opposite direction.

\item[(b)] \emph{Symmetry}: $\theta(x,y)=\theta(y,x)$ for every $x,y\in X$.

\item[(c)] \emph{Homogeneity}:
\begin{equation*}
\theta(ax,by)=\left\{
\begin{array}{ll}
\theta(x,y), & ab>0 \\
\pi -\theta(x,y),\text{ \ } & ab<0.%
\end{array}%
\right.
\end{equation*}

\item[(d)] \emph{Continuity}: If $x_{n}\rightarrow x$ and $y_{n}\rightarrow y$ (in the
norm), then $\theta(x_{n},y_{n})\rightarrow \theta(x,y)$.
\end{enumerate}

\medskip

\noindent In a normed space, the concept of angles between two vectors has been studied
intensively (see, for instance, \cite{Balestro, Diminnie2, Gunawan6, Horvath, Milicic2,
Thurev, Chen}). Here we shall be interested in the notion of angles between two subspaces
of a normed space using a semi-inner product.

\noindent Let $(X,\| \cdot \| )$ be a real normed space. The functional $g:X^{2}\rightarrow
\mathbb{R}$ defined by the formula
$$
g(x,y):=\frac{1}{2}\left\Vert x\right\Vert \left[ \tau _{+}(x,y)+\tau _{-}(x,y)\right],
$$
with
\begin{equation*}
\tau _{\pm }(x,y):=\lim_{t \rightarrow \pm 0}\frac{\|x+ty\| - \|x\| }{t},
\end{equation*}
satisfies the following properties:
\begin{enumerate}
\item[(1)] $g(x,x) =\| x\|^{2}$ for every $x\in X$;

\item[(2)] $g(ax,by)= ab \cdot g(x,y)$ for every $x,y\in X$ and $a,b\in \mathbb{R}$;

\item[(3)] $g(x,x+y)=\| x\|^{2}+g(x,y)$ for every $x,y\in X$;

\item[(4)] $| g(x,y)| \leq \|x\|\cdot\|y\|$ for every $x,y\in X$.
\end{enumerate}

\medskip

If, in addition, the functional $g(x,y)$ is linear in $y$, then $g$ is called
a \emph{semi-inner product} on $X$. For example, consider the space $\ell^{p}$
($1\leq p<\infty$) with the norm $\|x\|_{p}:=\left[ \sum_{k=1}^\infty
| \xi_{k}|^{p}\right] ^{\frac{1}{p}}$, $x=(\xi_k)$. Then the functional
\begin{equation*}
g(x,y):=\|x\|_{p}^{2-p}\sum_{k=1}^\infty |\xi_{k}|^{p-1}\text{sgn}\left( \xi_{k}\right)
\eta_{k}, \quad x:=\left(\xi _{k}\right),\ y:=\left(\eta _{k}\right) \in \ell^{p},
\end{equation*}
is a semi-inner product on $\ell^{p}\ (1\le p<\infty)$ \cite{Giles, Gunawan6}.
Note that, in general, $g$ is not commutative.

\bigskip

\noindent Using a semi-inner product $g$, Mili\`ci\`c \cite{Milicic1} introduced
the notion of \emph{$g$-orthogonality} on $X$, namely $x$ is said to be $g$-orthogonal
to $y$, denoted by  $x\perp _{g}y$, if $g(x,y)=0$. Note that in an inner product space,
the functional $g(x,y)$ is identical with the inner product
$\left\langle\cdot,\cdot\right\rangle$, and so the $g$-orthogonality coincides with
the usual orthogonality. In this article, we will develop the notion of $g$-angles
between two subspaces of a normed space and discuss its properties. We will begin our
discussion by studying the $g$-angle between two vectors in a normed space.

\section{Main Results}

\subsection{The $g$-angle between two vectors}

From now on, let $(X,\|\cdot\|)$ be a real normed space, unless otherwise stated.
In connection with the notion of $g$-ortogonality, we define the \emph{$g$-angle}
between two nonzero vectors $x$ and $y$ in $X$, denoted by $A_{g}(x,y)$, by the formula
\begin{equation*}
A_{g}(x,y):=\arccos \frac{g(y,x)}{\| x\| \cdot \left\Vert
	y\right\Vert }.
\end{equation*}
Note that $A_g(x,y) = \frac{1}{2}\pi$ if and only if $g(y,x) = 0$ or $y\perp _{g}x$.
If $X$ is an inner product space, the $g$-angle in $X$ is identical with the usual angle.

\bigskip

\begin{proposition}\label{proposition:2.1}
The $g$-angle $A_g(\cdot,\cdot)$ satisfies the following properties:

\begin{enumerate}
\item[\emph{(a)}]  If $x$ and $y$ are of the same direction, then  $A_{g}(x,y)=0$;
        if $x$ and $y$ are of opposite direction, then $A_{g}(x,y)=\pi $
        (part of the parallelism property).
		
\item[\emph{(b)}] $A_{g}(ax,by)=A_{g}(x,y)$ if $ab>0$; $A_{g}(ax,by)=\pi-A_{g}(x,y)$
        if $ab<0$ (the homogeneity property).
		
\item[\emph{(c)}] If $x_{n}\rightarrow x$ (in norm),
		then $A_{g}(x_{n},y)\rightarrow A_{g}(x,y)$ (part of the continuity property).
\end{enumerate}
\end{proposition}

\medskip

\noindent{\it Proof}.
	
\begin{enumerate}
\item[(a)] Let $y=kx$ for an arbitrary nonzero vector $x$ in $X$ and $k\in
		\mathbb{R}-\{0\}.$ We have%
		\begin{align*}
			A_{g}(x,y) = \arccos \frac{g(kx,x)}{\|x\| \cdot
				\| kx\|}=\arccos \frac{k\cdot g(x,x)}{|
				k| \| x\|^{2}} =\arccos \frac{k\|x\|^{2}}{|k|\|x\|^{2}}.
		\end{align*}%
		If $y=kx$ with $k>0$, then  $A_{g}(x,y)=\arccos (1)=0.$ If $y=kx$ with
        $k<0$, then $A_{g}(x,y)=\arccos (-1)=\pi.$
		
\item[(b)] Let  $a$ and $b\in\mathbb{R}-\{0\}.$ Observe that
		\begin{align*}
			A_{g}(ax,by)=\arccos \frac{ab\cdot g(y,x)}{|ab|( \|
				x\| \cdot \|y\|)}.
		\end{align*}%
		If $ab>0$, then $A_{g}(ax,by)=A_{g}(x,y).$ Likewise, if $ab<0$,
        then $A_{g}(ax,by)=\arccos \left( -\frac{g(y,x)}{
			\| x\|\cdot\|y\|}\right).$ Hence $A_{g}(ax,by)=\pi -A_{g}(x,y)$.
		
\item[(c)] If $x_{n}\rightarrow x$ (in norm), then
		\begin{align*}
		 |g(y,x_{n}-x)|\leq \|y\|\cdot\|x_{n}-x\|\longrightarrow 0.
		\end{align*}
		Observe that $g(y,x_{n}-x) =g(y,x_{n})-g(y,x)$. We have $g(y,x_{n})
        \longrightarrow g(y,x)$. Hence
		\begin{align*}
		A_{g}(x_{n},y)\rightarrow A_{g}(x,y),
		\end{align*}
as desired. \qed
\end{enumerate}

\bigskip

\remark Since $g$, in general, is not commutative, the $g$-angle does not satisfy
the symmetry property. For instance, in $\ell^{1}$ with $g(y,x):=\|y\|_{1}\sum_{k=1}^\infty
\text{sgn}(\eta _{k}) \xi_{k}$, take $x:=(1,1,0,\dots)$ and $y:=(-1,2,0,\dots)$,
so that we have $g(y,x)=0\neq  g(x,y)=2$. Likewise, the $g$-angle does not satisfy
the continuity property. For instance, in $\ell^{1}$ with the above functional $g$,
take $y_{n}:=(\frac{1}{n},1,0,\dots),\ x_{n}:=(1+\frac{1}{n},1,0, \dots),
y:=(0,1,0,\dots )$ and $x:=(1,1,0,\dots)$, so that we obtain $g(y_{n},x_{n})
\nrightarrow g(y,x)$.

\bigskip

\subsection{The $g$-angle between a $1$-dimensional subspace and a $t$-dimensional subspace}

Here, using a semi-inner product $g$, we will discuss the notion of $g$-angles between
two subspaces of a normed space. We first state the connection between the Gram determinant
$\Gamma(x_{1},\ldots ,x_{n}):=\det [g(x_{i},x_{k})]$, where $g(x_{i},x_{k})$ is the
$k$-th element of the $i$-th row, and the linearly independence of
$\left\{ {x_{1},\ldots, x_{n}}\right\} $ as in the following theorem.

\bigskip

\begin{theorem}\label{theorem:2.2}
Let $g$ be a semi-inner product in $X$ and
$\left\{ {x_{1},\dots, x_{n}}\right\}$ $\subset X$. If $\Gamma (x_{1},\dots, x_{n})
\neq 0,$ then $\left\{ {x_{1},\dots, x_{n}}\right\} $ is linearly independent.
\end{theorem}

\medskip

\begin{proof}
Suppose, on the contrary, that $\left\{{x_{1},\dots x_{n}}\right\}$ is linearly dependent.
Then there is an index $j$ with $1\leq j\leq n$ so that $x_{j}$ is a linear combination of
$x_{1},\dots,x_{j-1},x_{j+1},\dots,x_{n}$. Here the $j$-th column of $\Gamma$ is a linear
combination of the other columns. As a consequence, we have $\Gamma (x_{1},\dots,x_{n})=0$.
Hence, $\left\{{x_{1},\dots, x_{n}}\right\} $ must be a linearly independent set.
\end{proof}

\bigskip

\remark The converse of this theorem is not true. For example, take
$x_{1}:=(1,2,0,\dots)$ and $x_{2}:=(2,1,0,\dots)$ in $\ell^{1}$ with the usual
semi-inner product $g$. Clearly $x_{1}$ and $x_{2}$ are linearly independent.
But one may check that
$$g(x_{i},x_{j})=\|x_{i}\|_{1}\sum_{k=1}^\infty \text{sgn}(x_{ik})x_{jk}=9$$
for $i,j=1,2$, and hence $\Gamma (x_{1},x_{2})=\left| \begin{array}{cc}
9&9\\9&9\end{array}\right|=0$.

\bigskip

\noindent We shall now define the $g$-orthogonal projection of $y$ on subspace $S$ as follows.

\bigskip

\begin{definition}\label{fact:2.3}
\cite{Milicic0} Let $y$ be a vector of $X$ and $S=\text{span}\{x_{1},\dots,x_{n}\}$ be
subspace of $X$ with $\Gamma (x_{1},\dots, x_{n})=\det [g(x_{i},x_{k})]\neq 0$.
The $g$-{\it orthogonal projection} of $y$ on $S$, denoted by $y_{S}$, is defined by
\begin{align*}
y_{S}:=-\frac{1}{\Gamma (x_{1},\dots, x_{n})}\left\vert
\begin{array}{cccc}
0 & x_{1} & \cdots  & x_{n} \\
g(x_{1},y) & g(x_{1},x_{1}) & \cdots  & g(x_{1},x_{n}) \\
\vdots  & \vdots  & \ddots  & \vdots  \\
g(x_{n},y) & g(x_{n},x_{1}) & \cdots  & g(x_{n},x_{n})%
\end{array}%
\right\vert,
\end{align*}
and its $g$-{\it orthogonal complement} $y-y_{S}$ is given by
\begin{align*}
y-y_{S}=\frac{1}{\Gamma (x_{1},\dots, x_{n})}\left\vert
\begin{array}{cccc}
y & x_{1} & \cdots  & x_{n} \\
g(x_{1},y) & g(x_{1},x_{1}) & \cdots  & g(x_{1},x_{n}) \\
\vdots  & \vdots  & \ddots  & \vdots  \\
g(x_{n},y) & g(x_{n},x_{1}) & \cdots  & g(x_{n},x_{n})%
\end{array}%
\right\vert.
\end{align*}
\end{definition}
Note that the notation of the determinant $|\cdot|$ here has a special meaning
because the elements of the matrix are not in the same field. Since
\begin{align*}
g(x_{i},y-y_{S}) =\frac{1}{\Gamma( x_{1},\dots, x_{n}) }%
\left\vert
\begin{array}{cccc}
g\left( x_{i},y\right) & g\left( x_{i},x_{1}\right) & \cdots & g\left(
x_{i},x_{n}\right) \\
g\left( x_{1},y\right) & g\left( x_{1},x_{1}\right) & \cdots & g\left(
x_{1},x_{n}\right) \\
\vdots & \vdots & \ddots & \vdots \\
g\left( x_{n},y\right) & g\left( x_{n},x_{1}\right) & \cdots & g\left(
x_{n},x_{n}\right)%
\end{array}%
\right\vert,
\end{align*}%
we obtain $x_{i}\perp _{g}y-y_{S}$ for every $i=1,\dots,n$. For example, if
$S=\text{span}\left\{ x\right\}$,
then the $g$-orthogonal projection of $y$ on $x$ is
\begin{equation*}
y_{x}=\frac{g(x,y)}{\|x\|^{2}}x,
\end{equation*}%
and $y-y_{x}$ is the $g$-orthogonal complement $y$ on $x$. Notice here that
$x\perp _{g}y-y_{x}$. Note that the $g$-angle between two vectors in a normed
space is also the $g$-angle between these two vectors in the subspace spanned
by them.

Next, let $x_{1},\dots,x_{n}\in X$ be a finite sequence of linearly independent
vectors. We can construct a {\it left $g$-orthonormal sequence} $x_{1}^{\ast
},\dots,x_{n}^{\ast }$ with $x_{1}^{\ast }:=\frac{x_{1}}{\|x_{1}\|}$ and
\begin{equation}
x_{k}^{\ast }:=\frac{x_{k}-\left( x_{k}\right)_{s_{k-1}}}
{\|x_{k}-\left( x_{k}\right)_{s_{k-1}}\|},
\end{equation}%
where $S_{k-1}=\text{span}\left\{ x_{1}^{\ast },\dots,x_{k-1}^{\ast }\right\}$,
$k=2,\dots,n$.
We observe that $x_{k}^{\ast}\perp _{g}x_{l}^{\ast}$ for $k,l=1,\dots,n$ with $k<l$)
(see \cite{Gunawan7,Milicic0}).

Using the $g$-orthogonal projection, we define the $g$-angle between a
$1$-dimensional subspace $U =\text{span}\{u\}$ and a $t$-dimensional subspace
$V =\text{ span}\{v_{1},\dots, v_{t}\}$ of $X$ with $\Gamma (v_{1},\dots ,v_{t})\neq 0$
and $t\geq 1$ by
\begin{equation}\label{2}
\cos ^{2} A_{g}(U,V) :=\frac{(g(u_{V},u)) ^{2}}{\|u\|^{2}\|u_{V}\|^{2}},
\end{equation}
where $u_{V}$ denote the $g$-ortogonal projection of $u$ on $V$. Note that if
$U\subseteq V$, then $A_g(U,V)=0$. One may observe that if $X$ is an
inner product space, then
the definition of the $g$-angle in (\ref{2}) and the usual definition of angle
between two subspaces of $X$ are equivalent.

If we now write $u=u_{V}+u_{V}^{\perp }$ with $u_{V}^{\perp }$
is the $g$-orthogonal complement of $u$ on $V$, then $(2)$ becomes
\begin{equation*}
	\cos ^{2} A_{g}(U,V) =\frac{\|u_{V}\|^{2}}{\|u\|^{2}},
\end{equation*}
\noindent which tells us that the value of $\cos A_{g}(U,V)$ is equal to the ratio between
the `length' of the $g$-orthogonal projection of $u$ on $V$ and the `length' of $u$.
If $X=\ell^{p}$ with the semi-inner product $g$, then
$$
\cos ^{2} A_{g}(U,V)=\frac{\|u_{V}\|_{p}^{2}}{\|u\|_{p}^{2}}.
$$
Therefore, an explicit formula for the cosine of the $g$-angle between a $1$-dimensional
subspace $U=\text{span}\{u\}$ and $t$-dimensional subspace $V=\text{span}\{v_{1},\dots,v_{t}\}$
of $\ell^{p}$ can be presented as follows.

\bigskip

\begin{fact}\label{fact:2.5}
	If $U=\text{span}\{u\}$ is a $1$-dimensional subspace and $V=\text{span}
	\{v_{1},\dots ,v_{t}\}$ is a $t$-dimensional subspace of $\ell ^{p}$ with
   $\Gamma(v_{1},\dots,v_{t})\neq 0$, then
\begin{small}
		\begin{align*}
		\cos ^{2}A_{g}(U,V)=\left[ \underset{j_{t+1}}{\sum }\left\vert
		\sum\limits_{j_{t}}\cdots \sum\limits_{j_{1}}\left( \frac{1}{\Vert u\Vert
			_{p}}\prod\limits_{i=1}^{t}|v_{ij_{i}}^{\ast }|^{p-1}\text{sgn}%
		(v_{ij_{i}}^{\ast })\right) \left\vert
		\begin{array}{cccc}
		v_{1j_{1}} & \cdots  & v_{1j_{t}} & v_{1j_{t+1}} \\
		\vdots  & \ddots  & \vdots  & \vdots  \\
		v_{tj_{1}} & \cdots  & v_{tj_{t}} & v_{tj_{t+1}} \\
		u_{j_{1}} & \cdots  & u_{j_{t}} & 0%
		\end{array}%
		\right\vert \right\vert ^{p}\right] ^{\frac{2}{p}},
		\end{align*}%
	\end{small}
where in each summation the index ranges from 1 to $\infty$.
\end{fact}

\medskip

\noindent{\it Proof.} Suppose that $V=\text{span}\left\{v_{1},\dots,v_{t}\right\} $ with
$\Gamma(v_{1},\dots,v_{t})\neq 0.$ According to Theorem \ref{theorem:2.2},
$\{v_{1},\dots,v_{t}\}$ is linearly independent. Using $v_{1}^{\ast }=\frac{v_{1}}{\|v_{1}\|}$
 and so forth as in (1), we obtain  left
$g$-orthonormal set $\left\{ v_{1}^{\ast },\dots,v_{t}^{\ast }\right\}$. Notice that
$\text{span}\{ v_{1},\ldots,v_{t}\}=\text{span}\{v_{1}^{\ast },\ldots,v_{t}^{\ast }\}$. Hence
\begin{align*}
	u_{V}=-\frac{1}{\Gamma \left( v_{1}^{\ast },\dots,v_{t}^{\ast }\right) }
	\left\vert
	\begin{array}{cccc}
	0 & v_{1}^{\ast } & \cdots & v_{t}^{\ast } \\
	g(v_{1}^{\ast },u) & g(v_{1}^{\ast },v_{1}^{\ast })
	& \cdots & g(v_{1}^{\ast },v_{t}^{\ast }) \\
	\vdots & \vdots & \ddots & \vdots \\
	g(v_{t}^{\ast },u) & g(v_{t}^{\ast },v_{1}^{\ast })
	& \cdots & g(v_{t}^{\ast },v_{t}^{\ast })%
	\end{array}
	\right\vert.
\end{align*}

\noindent Observe that $\Gamma(v_{1}^{\ast},\dots,v_{t}^{\ast})=1$,
and so
\begin{align*}
\|u_{V}\|_{p}=& \,\left( \underset{j_{t+1}}{\sum }\left\vert\,
\left\vert
\begin{array}{cccc}
0 & v_{1j_{t+1}}^{\ast } & \cdots & v_{tj_{t+1}}^{\ast } \\
g(v_{1}^{\ast },u) & g(v_{1}^{\ast },v_{1}^{\ast }) & \cdots  &
g(v_{1}^{\ast },v_{t}^{\ast }) \\
\vdots  & \vdots  & \ddots  & \vdots  \\
g(v_{t}^{\ast },u) & g(v_{t}^{\ast },v_{1}^{\ast }) & \cdots  &
g(v_{t}^{\ast },v_{t}^{\ast })%
\end{array}%
\right\vert\, \right\vert ^{p}\right) ^{\frac{1}{p}} \\
=& \,\left( \underset{j_{t+1}}{\sum }\left\vert\, \left\vert
\begin{array}{cccc}
g(v_{1}^{\ast },v_{1}^{\ast }) & \cdots  & g\left( v_{1}^{\ast},v_{t}^{\ast}
\right)  & g(v_{1}^{\ast },u) \\
\vdots  & \ddots  & \vdots  & \vdots  \\
g(v_{t}^{\ast },v_{1}^{\ast }) & \cdots  & g(v_{t}^{\ast},v_{t}^{\ast}) &
g(v_{t}^{\ast },u) \\
v_{1j_{t+1}}^{\ast } & \cdots  & v_{tj_{t+1}}^{\ast } & 0%
\end{array}%
\right\vert\, \right\vert ^{p}\right) ^{\frac{1}{p}} \\
=& \,\left( \underset{j_{t+1}}{\sum }\left\vert\, \left\vert
\begin{array}{cccc}
g(v_{1}^{\ast },v_{1}^{\ast }) & \cdots  & g(v_{t}^{\ast },v_{1}^{\ast }) &
v_{1j_{t+1}}^{\ast } \\
\vdots  & \ddots  & \vdots  & \vdots  \\
g(v_{1}^{\ast },v_{t}^{\ast }) & \cdots  & g(v_{t}^{\ast },v_{t}^{\ast }) &
v_{tj_{t+1}}^{\ast } \\
g(v_{1}^{\ast },u) & \cdots  & g(v_{t}^{\ast },u) & 0%
\end{array}%
\right\vert\, \right\vert ^{p}\right) ^{\frac{1}{p}}.
\end{align*}

\noindent Hence, we have
	\begin{align*}
	\cos ^{2}A_{g}(U,V)=\frac{\|u_{V}\|_{p}^{2}}{\|
		u\|_{p}^{2}}=\left( \underset{j_{t+1}}{\sum }\left\vert\frac{1}{%
		\|u\|_{p}}\left\vert
	\begin{array}{cccc}
	g\left( v_{1}^{\ast },v_{1}^{\ast }\right)  & \cdots  & g\left( v_{t}^{\ast
	},v_{1}^{\ast }\right)  & v_{1j_{t+1}}^{\ast } \\
	\vdots  & \ddots  & \vdots  & \vdots  \\
	g(v_{1}^{\ast },v_{t}^{\ast })  & \cdots  & g(v_{t}^{\ast
	},v_{t}^{\ast })  & v_{tj_{t+1}}^{\ast } \\
	g(v_{1}^{\ast },u)  & \cdots  & g(v_{t}^{\ast },u)
	& 0%
	\end{array}%
	\right\vert \right\vert^{p}\right) ^{\frac{2}{p}}.
	\end{align*}
	
\noindent Since $v_{1}^{\ast }=v_{1}$ and $g(x,y)$ is linear in $y$, we have
	\begin{align*}
	\cos ^{2}A_{g}(U,V)=\left( \underset{j_{t+1}}{\sum }\left\vert \frac{1}{\|u\|_{p}}
    \left\vert
	\begin{array}{cccc}
	g(v_{1}^{\ast },v_{1})  & \cdots  & g(v_{t}^{\ast
	},v_{1})  & v_{1j_{t+1}} \\
	\vdots  & \ddots  & \vdots  & \vdots  \\
	g(v_{1}^{\ast },v_{t})  & \cdots  & g(v_{t}^{\ast
	},v_{t})  & v_{tj_{t+1}} \\
	g(v_{1}^{\ast },u)  & \cdots  & g(v_{t}^{\ast },u)
	& 0%
	\end{array}%
	\right\vert \right\vert^{p}\right)^{\frac{2}{p}}.
	\end{align*}%
Next, substituting $g(v_{i}^{\ast },v_{k}) =\|v_{i}^{\ast}\|_{p}^{2-p}
\sum\limits_{j_{i}}|v_{ij_{i}}^{\ast }|^{p-1}\text{sgn}(v_{ij_{i}}^{\ast })v_{kj_{i}}$
and using properties of determinants, we obtain
\begin{small}
	\begin{align*}
	\cos ^{2}A_{g}(U,V)=\left[ \underset{j_{t+1}}{\sum }\left\vert
	\sum\limits_{j_{t}}\cdots \sum\limits_{j_{1}}\left( \frac{1}{\Vert u\Vert
		_{p}}\prod\limits_{i=1}^{t}|v_{ij_{i}}^{\ast }|^{p-1}\text{sgn}%
	(v_{ij_{i}}^{\ast })\right) \left\vert
	\begin{array}{cccc}
	v_{1j_{1}} & \cdots  & v_{1j_{t}} & v_{1j_{t+1}} \\
	\vdots  & \ddots  & \vdots  & \vdots  \\
	v_{tj_{1}} & \cdots  & v_{tj_{t}} & v_{tj_{t+1}} \\
	u_{j_{1}} & \cdots  & u_{j_{t}} & 0%
	\end{array}%
	\right\vert \right\vert ^{p}\right] ^{\frac{2}{p}}.
	\end{align*}
\end{small}

\noindent This proves the fact. \qed

\bigskip

\example Consider $\ell^{1}$ with the usual semi-inner product $g$. Take $U=\text{span}%
\{u\}$ and $V=\text{span}\{v_{1},v_{2}\}$, with $u:=(1,2,1,0,\ldots),$
$v_{1}:=(1,0,0,0,\ldots)$, and $v_{2}:=(0,1,0,0,\ldots )$. We obtain
\begin{align*}
\cos ^{2}A_{g}(U,V)=& \,\frac{1}{16}\left[ \underset{j_{3}}{\sum }\left\vert
\sum\limits_{j_{2}}\sum\limits_{j_{1}}\text{sgn}(v_{1j_{1}})\text{sgn}%
(v_{2j_{2}})\left\vert
\begin{array}{ccc}
v_{1j_{1}} & v_{1j_{2}} & v_{1j_{3}} \\
v_{2j_{1}} & v_{2j_{2}} & v_{2j_{3}} \\
u_{j_{1}} & u_{j_{2}} & 0%
\end{array}%
\right\vert \right\vert \right] ^{2} \\
=& \,\frac{1}{16}\left[ \underset{j_{3}}{\sum }\left\vert \sum\limits_{j_{2}}%
\text{sgn}(v_{2j_{2}})\left\vert
\begin{array}{ccc}
1 & v_{1j_{2}} & v_{1j_{3}} \\
0 & v_{2j_{2}} & v_{2j_{3}} \\
1 & u_{j_{2}} & 0%
\end{array}%
\right\vert \right\vert \right] ^{2} \\
=& \,\frac{1}{16}\left[ \left\vert\, \left\vert
\begin{array}{ccc}
1 & 0 & 1 \\
0 & 1 & 0 \\
1 & 2 & 0%
\end{array}%
\right\vert\, \right\vert +\left\vert\, \left\vert
\begin{array}{ccc}
1 & 0 & 0 \\
0 & 1 & 1 \\
1 & 2 & 0%
\end{array}%
\right\vert\, \right\vert \right] ^{2}=\frac{9}{16}.
\end{align*}%
Hence, the $g$-angle between the two subspaces $U$ and $V$ is $\arccos(\frac{3}{4})$.

\bigskip

\subsection{The $g$-angle between a $2$-dimensional subspace and a $t$-dimensional subspace}

In this section, we discuss the $g$-angle between subspaces $U$ and $V$ where $U$ is
a $2$-dimensional subspace and $V$ is a $t$-dimensional subspace with $t\geq 2$.
First, we define the function $\Lambda (\cdot,\cdot)$ on $X \times X$ by
\begin{equation}
 \Lambda (x,y):=\left\vert
\begin{array}{cc}
\left\vert g(x,x)\right\vert & \left\vert g(x,y)\right\vert  \\
\left\vert g(y,x)\right\vert & \left\vert g(y,y)\right\vert%
\end{array}%
\right\vert^\frac{1}{2}. \label{1}
\end{equation}
Note that in a real inner product space $(X,\langle\cdot,\cdot\rangle)$,
$\Lambda (\cdot,\cdot)$ is identical with
the {\it standard $2$-norm} $\|\cdot,\cdot\|$ which is given by
\begin{equation*}
\|x,y\| :=\left\vert
\begin{array}{cc}
\left\langle x,x\right\rangle & \left\langle x,y\right\rangle \\
\left\langle y,x\right\rangle & \left\langle y,y\right\rangle%
\end{array}%
\right\vert ^{\frac{1}{2}}.
\end{equation*}%

The following proposition lists some properties of the function $\Lambda(\cdot,\cdot)$.

\bigskip

\begin{proposition}\label{proposition:2.6}
	The function $\Lambda(\cdot,\cdot)$ defined by (\ref{1}) satisfies the following properties:
\begin{enumerate}
	\item[\emph{(a)}] $\Lambda(x,y)\ge 0$  for every $x,y\in X$;
If $x$ and $y$ are linearly dependent, then $\Lambda(x,y)=0$;
	
	\item[\emph{(b)}] $\Lambda(x,y)=\Lambda(y,x)$ for every $x,y\in X$;
	
	\item[\emph{(c)}] $\Lambda(\alpha x,y) =\alpha \Lambda(x,y)$ for every $x,y\in X$
    and $\alpha\in \mathbb{R}$;
	
	\item[\emph{(d)}] $\Lambda(x,y)\le \| x\|\cdot\|y\|$ for every $x,y\in X$.
\end{enumerate}	
\end{proposition}

\medskip

\noindent{\it Proof.}
\begin{enumerate}
\item[(a)] Using properties of determinants and the functional $g$, we have
$$
	(\Lambda (x,y))^{2} =\,\|x\|^{2}\|
	y\|^{2}-|g(x,y)||g(y,x)|
	\ge \,\|x\|^{2}\|y\|^{2}-\|x\|^{2}\|y\|^{2}=0.
$$
Let $y=kx$ with $k\in\mathbb{R}$. Then we have
$\Lambda(ky,y)=\sqrt{\|ky\|^{2}\|y\|^{2}-|g(ky,y)||g(y,ky)|} =0.$
    	
\item[(b)] Observe that
\begin{align*}
	\Lambda (x,y)=&\,\sqrt{\|x\|^{2}\|y\|^{2}-|g(x,y)| |g(y,x)|} \\
	=&\,\sqrt{\|y\|^{2}\|x\|^{2}-|g(y,x)|| g(x,y)|} \\
	=&\,\Lambda (y,x).
\end{align*}

\item[(c)] Observe that
\begin{align*}
	\Lambda(\alpha x,y) =\sqrt{\|\alpha x\|^{2}\|y\|^{2}-|g(\alpha x,y)| |g(y,\alpha x)| }
=|\alpha| \Lambda (x,y).
\end{align*}
\item[(d)] Observe that
\begin{align*}
	\Lambda (x,y)=\sqrt{\|x\|^{2}\|y\|^{2}-|g(x,y)||g(y,x)| }\le \|x\|\cdot\| y\|,
\end{align*}
as desired.\qed
\end{enumerate}

\bigskip

\remark The converse of part (a) is not true. For instance, take $x:=(1,2,0,\dots)$ and
$y:=(2,1,0,\dots)$ in $\ell^{1}$ (with the semi-inner product $g$, as usual). Clearly $x$ and $y$
are linearly independent. But one may check that $\Lambda(x,y)=0$. Likewise, $\Lambda$
does not satisfy the triangle inequality. For example, take $x:=(3,1,0,\dots),\ y:=(-2,0,0,\dots),
\ z:=(0,2,0,\dots)$ in $\ell^{1}$. Then we have $\| x\| =4,$ $\|y\| =2,$ $\|z\| =2$, $\|y+z\|
=4,$ $g(x,y)=-8,$ $g(y,z)=-6$, $g(x,z)=8$, $g(z,x)=2$, $g(x,y+z)=0$, and $g(y+z,x)=-8.$
Hence $\Lambda(x,y)=4$, $\Lambda(x,z)=4\sqrt{3}$, and $\Lambda(x,y+z)=16$,
so that $\Lambda (x,y+z)>\Lambda (x,y)+\Lambda (x,z)$.

\bigskip

Using the function $\Lambda(\cdot,\cdot)$, we now define the $g$-angle between a $2$-dimensional
subspace $U :=\text{span}\{u_{1},u_{2}\}$ of $X$ with $\Lambda(u_{1},u_{2})\neq 0$
and a $t$-dimensional subspace $V :=\text{ span}\{v_{1},\dots, v_{t}\}$ of $X$
with $\Gamma(v_{1},\dots,v_{t})\neq 0$ ($t\geq 2$) by
\begin{equation}
\cos^{2}A_{g}(U,V):=\frac{(\Lambda (u_{1V},u_{2V}))^{2}}{%
	(\Lambda(u_{1},u_{2}))^{2}}   \label{4}
\end{equation}
where $u_{iV}$ denote the $g$-orthogonal projection of $u_{i}$'s on $V$ with $i=1,2$.
Note that in a standard $2$-normed space, the definition of $g$-angle in (\ref{4}) is
identical with the angle defined in \cite{Gunawan7}.

\medskip

According to the following proposition, the definition of $g$-angle in (\ref{4}) makes sense.

\bigskip

\begin{proposition}\label{proposition:2.7}
The ratio on the right hand side of (\ref{4}) is a number in $[0, 1]$, and it is independent
of the choice of basis for $U$ and $V$ .
\end{proposition}

\medskip

\begin{proof}
Assuming particularly that $\{u_{1},u_{2}\} $ is left orthonormal, we have
$\Lambda (u_{1},u_{2}) =1$ and
\begin{equation*}
	(\Lambda (u_{1V},u_{2V}) )^{2} =\|u_{1V}\|
	^{2}\|u_{2V}\|^{2}-|g(u_{1V},u_{2V})||g(u_{2V},u_{1V}|.
\end{equation*}
According to \cite{Milicic0}, we have $\|u_{iV}\|\le \|u_{i}\|$ for $i=1,2$.
Hence $(\Lambda(u_{1V},u_{2V}))^{2}\le 1$. Therefore, the ratio is a number in $[0,1]$.

\noindent Secondly, note that the $g$-orthogonal projection of $u_{i}$'s on $V$ is independent
of the choice of basis for $V$ \cite{Milicic0}. Moreover, since  $g$-orthogonal projections
are linear transformations, the ratio of (\ref{4}) is also invariant under any change of basis
for $U$. Indeed, the ratio is unchanged if swap $u_{1}$ and $u_{2}$, replace $u_{1}$ with
$u_{1}+\alpha u_{2}$, or replace $u_{1}$ with $\alpha u_{2}$ where $\alpha\neq 0$.
\end{proof}

\bigskip

\section{Concluding Remarks}

The formula (\ref{4}) can be used to compute the $g$-angle between two subspaces of $\ell^p$ as
follows. Let $V=\text{span}\{ v_{1},\dots,v_{t}\} $ with $\Gamma(v_{1},\dots,v_{t}) \neq 0.$
According to Theorem \ref{theorem:2.2}, $\{v_{1},\dots,v_{t}\} $ is linearly independent.
Using $v_{1}^{\ast }=\frac{v_{1}}{\|v_{1}\|}$ and so forth as in (1), we obtain the left
$g$-orthonormal set $\{v_{1}^{\ast },\dots,v_{t}^{\ast }\}$. Here
$\text{span}~\{v_1,\ldots,v_t\}=\text{span}~\{v_{1}^{\ast },\dots,v_{t}^{\ast }\}$.
Hence, for $i=1,2$, we have
$$
	u_{iV} =\,-\left\vert
	\begin{array}{cccc}
		0 & v_{1}^{\ast } & \cdots  & v_{t}^{\ast } \\
		g(v_{1}^{\ast },u_{i})  & g(v_{1}^{\ast },v_{1}^{\ast
		})  & \cdots  & g(v_{1}^{\ast },v_{t}^{\ast })  \\
		\vdots  & \vdots  & \ddots  & \vdots  \\
		g(v_{t}^{\ast },u_{i})  & g(v_{t}^{\ast },v_{1}^{\ast
		})  & \cdots  & g(v_{t}^{\ast },v_{t}^{\ast })
	\end{array}
	\right\vert
	=\,-\left\vert
	\begin{array}{cccc}
		g(v_{1}^{\ast },v_{1}^{\ast })  & \cdots  & g(v_{t}^{\ast
		},v_{1}^{\ast })  & v_{1}^{\ast } \\
		\vdots  & \ddots  & \vdots  & \vdots  \\
		g(v_{1}^{\ast },v_{t}^{\ast })  & \cdots  & g(v_{t}^{\ast
		},v_{t}^{\ast })  & v_{t}^{\ast } \\
		g(v_{1}^{\ast },u_{i})  & \cdots  & g(v_{t}^{\ast
		},u_{i})  & 0%
	\end{array}
	\right\vert.
$$
Since $v_{1}^{\ast }=v_{1}$ and $g(x,y)$ is linear in $y$, we obtain
\begin{align*}
u_{iV}=-\left\vert
\begin{array}{cccc}
g(v_{1}^{\ast },v_{1})  & \cdots  & g(v_{t}^{\ast},v_{1})  & v_{1} \\
\vdots  & \ddots  & \vdots  & \vdots  \\
g(v_{1}^{\ast },v_{t})  & \cdots  & g(v_{t}^{\ast},v_{t})  & v_{t} \\
g(v_{1}^{\ast },u_{i})  & \cdots  & g(v_{t}^{\ast},u_{i})  & 0%
\end{array}%
\right\vert.
\end{align*}%
Substituting $g(v_{k}^{\ast },v_{i}) =\| v_{k}^{\ast
}\|_{p}^{2-p}\sum\limits_{j_{k}}|v_{kj_{k}}^{\ast}|
^{p-1}\text{sgn}(v_{kj_{k}}^{\ast })v_{ij_{k}}$, we get
\begin{align*}
u_{iV}=-\sum\limits_{j_{t}}\cdots \sum\limits_{j_{1}}|
v^{*}_{1j_{t}}|^{p-1}\text{sgn}(v^{*}_{1j_{t}})\cdots |
v^{*}_{tj_{1}}|^{p-1}\text{sgn}(v^{*}_{tj_{1}})\left\vert
\begin{array}{cccc}
v_{1j_{1}} & \cdots  & v_{1j_{t}} & v_{1} \\
\vdots  & \ddots  & \vdots  & \vdots  \\
v_{tj_{1}} & \cdots  & v_{tj_{t}} & v_{t} \\
u_{ij_{1}} & \cdots  & u_{ij_{t}} & 0%
\end{array}%
\right\vert.
\end{align*}
Using this formula, the value of $\cos^{2}A_{g}(U,V)$ (and hence $A_g(U,V)$) can
be computed. For instance, in $\ell ^{1}$ (with the usual semi-inner product $g$),
let $U=\text{span}\{u_{1},u_{2}\}$ and $V=\text{span}\{v_{1},v_{2},v_{3}\}$
with $u_{1}:=(1,1,2,3,0,\ldots) $, $u_{2}:=(2,1,-3,2,0,\ldots) $, $v_{1}:=(1,0,0,0,0,\ldots )$,
$v_{2}:=(0,1,0,0,0,\ldots )$, and $v_{3}:=(0,0,1,0,0,\ldots )$. We obtain
$u_{1V}=\left( 1,1,2,0,0,\ldots \right)$ and $u_{2V}=(2,1,-3,0,0,\ldots)$.
Moreover, $\|u_{1}\| =7$, $\|u_{2}\|=8$, $g( u_{1},u_{2}) =14$, $g(u_{2},u_{1}) =24$,
$\|u_{1V}\|=4$, $\|u_{2V}\|=6$, $g(u_{1V},u_{2V}) =0$, and $g(u_{2V},u_{1V}) =0$.
Thus $\cos ^{2}A_{g}(U,V)=\frac{36}{167}$, so that
$A_{g}(U,V)=\arccos(\frac{6}{167}\sqrt{167})$.

\bigskip

\textbf{Acknowledgement}. The research is supported by ITB Research and
Innovation Program No. 107r/I1.C01/PL/2017. The authors thank the referee
for his/her useful comments and suggestions on the earlier version of
this paper.

\end{document}